\documentclass[12pt]{article}
\usepackage{amsthm}
\usepackage{amsmath}
\usepackage{amsfonts}
\usepackage{graphicx}
\usepackage{latexsym}
\usepackage{amssymb}

\DeclareMathAlphabet{\mathbfsl}{OT1}{ppl}{b}{it} 

\newcommand{\deff}{\mbox{$\stackrel{\rm def}{=}$}}

\makeatletter
 \DeclareRobustCommand{\nsbinom}{\genfrac[]\z@{}}
 \makeatother
 
 \newcommand{\sbinomq}[2]{\nsbinom{{#1}}{{#2}}_{q}}

\newcommand{\field}[1]{\mathbb{#1}}
\newcommand{\A}{\field{A}}
\newcommand{\dB}{\field{B}}

\newcommand{\F}{\field{F}}

\newcommand{\dS}{\field{S}}
\newcommand{\T}{\field{T}}

\newcommand{\cN}{{\cal N}}

\newcommand{\cG}{{\cal G}}

\newcommand{\linadd}{\kern1pt\mbox{\small$\boxplus$}\kern1pt}

\newtheorem{theorem}{Theorem}
\newtheorem{lemma}{Lemma}

\newtheorem{cor}{Corollary}

\begin{document}

\bibliographystyle{plain}

\title{
\begin{center}
On the Structure of the $q$-Fano Plane
\end{center}
}
\author{
{\sc Tuvi Etzion}\thanks{Department of Computer Science, Technion,
Haifa 32000, Israel, e-mail: {\tt etzion@cs.technion.ac.il}.}}

\maketitle

\begin{abstract}
The $q$-Fano plane is the $q$-analog of the Steiner system $S(2,3,7)$.
For any given prime power $q$ it is not known whether the $q$-Fano plane exists or not.
We consider the structure of the $q$-Fano plane for any given $q$ and
conclude that most of its structure is known. Even so, we were unable
to determine whether it exists or not. A special attention is given for the
case $q=2$ which was considered by most researchers before.
\end{abstract}

\vspace{0.5cm}


\vspace{0.5cm}


\footnotetext[1] { This research was supported in part by the Israeli
Science Foundation (ISF), Jerusalem, Israel, under
Grant 10/12.}

\newpage
\section{Introduction}

Let $\F_q$ be the finite field with $q$ elements
and let $\F_q^n$ be the set of all vectors
of length $n$ over~$\F_q$. $\F_q^n$ is a vector
space with dimension $n$ over $\F_q$. For a given integer $k$,
$1 \leq k \leq n$, let $\cG_q(n,k)$ denote the set of all
$k$-dimensional subspaces of $\F_q^n$. $\cG_q(n,k)$ is often
referred to as Grassmannian. It is well known that
$$ \begin{small}
| \cG_q (n,k) | = \sbinomq{n}{k}
\deff \frac{(q^n-1)(q^{n-1}-1) \cdots
(q^{n-k+1}-1)}{(q^k-1)(q^{k-1}-1) \cdots (q-1)}
\end{small}
$$
where $\sbinomq{n}{k}$ is the $q-$\emph{ary Gaussian
coefficient}~\cite[pp. 325-332]{vLWi92}.

A \emph{Steiner system} $S(t,k,n)$ is a collection $S$ of $k$-subsets
from an $n$-set $\cN$ such that each $t$-subset of $\cN$ is contained
in exactly one element of $S$. Steiner systems were subject to an
extensive research in combinatorial designs~\cite{CoDi07}.
A \emph{$q$-Steiner system} $\dS_q(t,k,n)$ is a collection $\dS$ of
elements from $\cG_q(n,k)$ (called \emph{blocks}) such that each element from
$\cG_q(n,t)$ is contained in exactly one block of $\dS$. There are some well-known~\cite{ScEt02,Suz90}
divisibility necessary conditions for the existence of the $q$-Steiner
system $\dS_q(t,k,n)$ which are implied by the following results.

\begin{theorem}
\label{thm:derived}
If a $q$-Steiner system $\dS_q (t,k,n)$ exists then
for each $i$, $1 \leq i \leq t-1$, a $q$-Steiner system
$\dS_q(t-i,k-i,n-i)$ exists.
\end{theorem}

\begin{cor}
If a $q$-Steiner system $\dS_q(t,k,n)$ exists then for all $0 \leq i \leq t-1$,
$$
\frac{\sbinomq{n-i}{t-i}}{\sbinomq{k-i}{t-i}}
$$
must be integers.
\end{cor}

While lot of information is known about the existence of Steiner systems~\cite{CoDi07,Kee14}, our knowledge of
$q$-Steiner systems is quite limited.
Until recently, the only known $q$-Steiner systems $\dS_q(t,k,n)$ were either trivial or
for $t=1$, where such structure exists if and only if $k$ divides $n$.
These systems are known as \emph{spreads} in projective geometry.
Recently, the first $q$-Steiner system $\dS_q(t,k,n)$
with $t \geq 2$ was found.
This is a $q$-Steiner system $\dS_2(2,3,13)$ which has a large
automorphism group~\cite{BEOVW}. Using derived and residual designs
it was proved that sometimes the necessary conditions for the existence
of a $q$-Steiner system $\dS_q(t,k,n)$ are not sufficient~\cite{KiLa15}.
The first set of parameters for which the existence question of $q$-Steiner systems
is not settled is the parameters for the $q$-analog
of the Fano plane, i.e. the $q$-Steiner systems $\dS_q(2,3,7)$, which
will be called also in this paper the $q$-Fano planes. It was proved
in~\cite{BKN15} that if such system exists for $q=2$ then it has an automorphism of small order.
A new approach to examine $q$-Steiner systems was suggested recently in our paper~\cite{Etz15}.
This approach is based on puncturing coordinates in the system and deriving sets of equations
that must be satisfied if such a system exists.

In this paper we examine the structure of the $q$-Fano plane based
on the puncturing method. We examine the puncturing
of the $q$-Steiner system $\dS_q(2,3,7)$ on different sets of coordinates
and combine the puncturing on these different sets of coordinates to obtain more
tight results on the structure of the $q$-Fano plane if exists. We will assume
throughout our discussion that the
reader is aware of the definitions and results presented in~\cite{Etz15},
even so a brief summary will be given.
The rest of this paper is organized as follows. In Section~\ref{sec:structure}
we prove our main results concerning the structure of the $q$-Fano plane.
The section is divided into three parts. In Subsection~\ref{sec:punctured}
we give a short introduction for punctured systems and discuss how it can help
to determine if the $q$-Fano plane exists or not. In Subsection~\ref{sec:q=2}
we consider the $q$-Steiner system $\dS_2(2,3,7)$ and discuss its structure.
In Subsection~\ref{sec:q>1} we discuss the structure of the $q$-Fano plane
for any prime power $q$. Although the discussion for $q=2$ does not provide much more
information than the one for general $q$, we feel it is better understood if first the
case $q=2$ is discussed. In Section~\ref{sec:2punctured} we consider one possible
structure of the 2-punctured $q$-Fano plane. The given construction is only one
of a few possibilities to construct this design. In Section~\ref{sec:conclude} we conclude with a few directions
for future research on how to advance the knowledge on this problem and maybe settle
it for good.

\section{The structure of the $q$-Fano plane}
\label{sec:structure}

In this section we present the structure of the $q$-Fano plane (if exists)
based on its punctured designs as defined in~\cite{Etz15}. We start with a brief
introduction to punctured designs in Subsection~\ref{sec:punctured}. In Subsection~\ref{sec:q=2}
the structure of the $q$-Fano plane for the binary field, i.e. $q=2$ is introduced.
These results are straightforward generalized for any prime power $q$ in Subsection~\ref{sec:q>1}.
The $q$-Fano plane $\dS_2(2,3,7)$ is the one on which most research was done in the past, e.g.~\cite{BKN15,EtVa11a,Tho96}
and this is the first reason we discuss it first even so the generalization is straightforward.
The second reason is that its understanding might be slightly simpler and it will make
a good preparation for the discussion for any prime power $q$. The third reason is that its size
is smaller and hence with the mentioned figures for various substructures of
the $q$-Fano plane, one can take it as a toy example to
try and construct it by hand, needless to say it might be easier to check its existence with computer search.
Finally, note that sometimes we have to consider for $q>2$ 1-subspaces instead of vectors for $q=2$, in
this generalization. This is another reason for the separation that made between $q=2$ and general $q$.

\subsection{Punctured $q$-Steiner systems}
\label{sec:punctured}

Punctured codes is a well-known concept in coding theory.
Given an $n \times m$ array $A$, the \emph{punctured} array $A'$ is an $n \times (m-1)$
array obtained from $A$ by deleting one of the columns of $A$.
For a subspace $X \in \cG_q(n,k)$ the vectors in $X$ will be regarded as row
vectors and any coordinate in all the vectors will be viewed together as a column of the subspace.
For a subspace $X \in \cG_q(n,k)$ the punctured subspace $X'$ by the $i$th coordinate of $X$ is defined
as the subspace obtained by deleting coordinate $i$ from all the vectors of $X$.
If not specified, the deleted coordinate should
be understood from the context. If it won't be understood from the context then it is
assumed that the punctured coordinate is the last one.
The result of this puncturing is a new subspace of
$\F_q^{n-1}$. If $X$ does not contain the unity vector with a \emph{one} in the $i$th coordinate
then $X'$ is a subspace in $\cG_q(n-1,k)$. If $X$ contains the unity vector with a \emph{one} in the $i$th coordinate
then $X'$ is a subspace in $\cG_q(n-1,k-1)$.
Assume that we are given a set $\dS$ of $k$-subspaces from $\cG_q(n,k)$. The \emph{punctured set}
$\dS'$ is defined as $\dS' = \{ X' ~:~ X \in \dS \}$ and it can contain elements only from
$\cG_q(n-1,k)$ and $\cG_q(n-1,k-1)$. The set $\dS'$ is regarded as multiset and hence $|\dS'|=|\dS|$.
A subspace can be punctured several times. A $k$-subspace $X$, of $\F_q^n$, is punctured $p$ times to a $p$-\emph{punctured}
subspace $Y$ of $\F_q^{n-p}$. The subspace $Y$ can be an $s$-subspace for any $s$ such that $\max \{ 0, k-p \} \leq s \leq \min \{ k, n-p \}$.
Similarly, we define a $p$-\emph{punctured} set.

For puncturing there is an inverse operation called \emph{extension}.
A $t$-subspace $X$ of $\F_q^m$ is \emph{extended} to a $t'$-subspace $Y$ of $\F_q^{m'}$, where
$t' \geq t$, $m'>m$, and $m'-m \geq t'-t$, if $X$ is the subspace obtained from $Y$ by
puncturing $Y$ $m'-m$ times. Note, that the extension of a subspace is not always unique (see Lemma~\ref{lem:addNone}),
but there is always a unique outcome for puncturing. In other words, a
$p$-punctured $t$-subspace $X$ of $\F_q^m$
can be obtained from a few different $t'$-subspaces of $\F_q^{m+p}$.
The following results concerning extension were proved in~\cite{Etz15}.

\begin{lemma}[\bf from $t$-subspace to $t$-subspace, one extension]
\label{lem:addNone}
If $X$ is a $t$-subspace in $\F_q^n$ then it can be extended
in exactly $q^t$ distinct ways to a $t$-subspace of $\F_q^{n+1}$.
\end{lemma}

\begin{lemma}[\bf from $t$-subspace to $(t+1)$-subspace, one extension]
\label{lem:addOne}
If $X$ is a $t$-subspace of $\F_q^n$ then it can be extended
in exactly one way to a $(t+1)$-dimensional subspace of $\F_q^{n+1}$.
\end{lemma}

A $p$-punctured $q$-Steiner system $\dS_q(t,k=t+1,n;m)$, $m=n-p$, is a system $\dS$ of subspaces from $\F_q^m$,
in which each $t$-subspace of $\F_q^n$ can be obtained exactly once by extending
$p$ times all the subspaces of $\dS$, where the appearances of the same
subspace of $\F_q^m$ in $\dS$ are extended in parallel for this purpose.

A system of equations that must be satisfied based on the definition of
a $p$-punctured $q$-Steiner system can be obtained.
Each $s$-subspace $X$ of $\F_q^m$, $\max \{0,t-p\} \leq s \leq \min \{t,m\}$, yields
one equation related to the way it is covered in $\dS_q(t,k=t+1,n;m)$.
Each $r$-subspace $Y$ of $\F_q^m$, $\max \{ s, \max \{0,k-p\}\} \leq r \leq \min \{s+1,m\}$,
yields one variable which is the number of appearances of $Y$ in
the $p$-punctured $q$-Steiner system $\dS_q(t,k=t+1,n;m)$.
The set of equations relate to the way in which each punctured $t$-subspaces of $\F_q^m$ is
covered by the punctured $r$-subspaces.
The easiest designs to construct based on the system of equations are those
for which all the $r$-subspaces of $\F_q^m$, for a fixed $r$, appear with the same amount
of times in the design. For such a
type of a design $\T$ we define a set of $\min \{ m,k \} +1$ variables; $X_{0,m}$ is the number
of 0-subspaces in $\T$; $X_{1,m}$ is the number of times each 1-subspace of $\F_q^n$ appears in $\T$;
$X_{2,m}$ is the number of times each 2-subspace of $\F_q^n$ appears in $\T$; and so on. Such a design
will be called \emph{uniform}. If there is no solution to the general (not the uniform one) set of equations then also
the $q$-Steiner system $\dS_q(t,k=t+1,n)$ which was punctured does not exist. If there is a solution, it might throw some light
on the structure of the $q$-Steiner system $\dS_q(t,k=t+1,n)$.

\subsection{The Structure of $\dS_q(2,3,7)$, $q=2$}
\label{sec:q=2}

Throughout our discussion, let $\dS$ be a $q$-Steiner system $\dS_2(2,3,7)$.
We start with a uniform solution for the 3-punctured $q$-Steiner system
$\dS_2(2,3,7;4)$. Such a uniform solution, given in~\cite{Etz15} implies that $X_{0,4}=1$ which implies
that $X_{1,4}=0$, $X_{2,4}=4$, and $X_{3,4}=16$. We set $X_{0,4}=1$ w.l.o.g.,
since in any system one subspace can be chosen w.l.o.g.~. Furthermore, in this case
$X_{0,4}=1$ implies that the design is uniform~\cite{Etz15}. $X_{0,4}=1$ implies that the 3-subspace whose
first four columns are all-zero columns is contained in $\dS$. Let $Z_1$ denote this
3-subspace in which the first four columns are all-zero. For symmetry let $Z_2$ denote
the 3-subspace of $\F_2^7$ in which the last four columns are all-zero columns.

In $\dS$, each nonzero vector of $\F_2^7$ appears in exactly 21 3-subspaces.
By puncturing the last coordinate of each 3-subspace of $\dS$,
all the 3-subspaces which contain the vector 0000001 will be punctured
into a spread, i.e. a $q$-Steiner system $\dS_2(1,2,6)$.

Next, in our exposition we exclude $Z_1$ from $\dS$ for the current paragraph.
Each 3-subspace of $\dS$ which contains a nonzero vector which starts with four \emph{zeroes}
is 3-punctured into a 2-subspace of $\F_2^4$. There are 7 such vectors, each one appears
in twenty distinct 3-subspaces of $\dS \setminus \{ Z_1 \}$ (since the only 3-subspace of $\dS$ which
contains two of them is $Z_1$; thus, clearly two of them cannot appear together in the same 3-subspace of $\dS \setminus \{ Z_1 \}$)
for a total of 140 such 3-subspaces which are punctured into the 140 (non-distinct) 2-subspaces
of the 3-punctured $q$-Steiner system $\dS_2(2,3,7;4)$ derived from $\dS$.
Since each nonzero vector of $\F_2^7$ appears in exactly one 3-subspace with each
other nonzero vector of $\F_2^7$, it follows that each such vector is responsible
for exactly twenty 2-subspaces (some of them are identical) of the 3-punctured $q$-Steiner system $\dS_2(2,3,7;4)$.
Hence, each nonzero prefix of length four of vectors of $\F_2^7$ (there are fifteen such nonzero prefixes) appears $\frac{20 \cdot 3}{15}=4$
times in these twenty 2-subspaces. There are many possible partitions to obtain such seven 3-punctured sets
of twenty 2-subspaces for this purpose.
One suggestion was given in the recursive construction of~\cite{Etz15}.
It will be discussed again in Section~\ref{sec:2punctured}.

We continue by imposing w.l.o.g. a certain structure on $\dS$.
Many structures can be imposed, each one will imply a different specifications for $\dS$. The one we suggest
now seems to be one of the most successful ones. We already forced the 3-subspace $Z_1$
to be a 3-subspace in $\dS$. Now, we assume
that w.l.o.g. also the 3-subspace whose last four columns are all-zero, i.e. $Z_2$, is also a 3-subspace in $\dS$.
Next, we will show why it can be assumed that $Z_2 \in \dS$ without loss of generality.
In the 3-punctured $q$-Steiner system $\dS_2(2,3,7;4)$ each 3-subspace of $\F_2^4$ appears exactly 16 times.
There are fifteen such 3-subspaces of $\F_2^4$, for one of them, the first three columns
form the unique 3-subspace in $\F_2^3$, and these three columns can be followed by the all-zero column (which is one
of the eight possible extensions of the unique 3-subspace of $\F_2^3$ by Lemma~\ref{lem:addNone}). Let $X$
be a 3-subspace of $\dS$ which is 3-punctured to these 4 columns. To continue
our discussion, we need the following simple lemma

\begin{lemma}
\label{lem:linear_comb}
If $\hat{\dS}$ is a $q$-Steiner system $\dS_q(t,k,n)$ then the system obtained by replacing
the $j$th column (for any $j$, $1\leq j \leq n$) in all the $k$-subspaces of $\hat{\dS}$,
by a linear combination of the columns which contains the $j$th column (in any nonzero
multiplicity), is also a $q$-Steiner system $\dS_q(t,k,n)$.
\end{lemma}

Therefore, since the first three columns of $X$ have rank three it follows that we can
form some specific three linear combinations, containing the 5th, the 6th, and the 7th column of $X$, respectively.
Each such linear combination will sum
to \emph{zero} for the related column of $X$. We replace the 5th, 6th, and 7th columns of $X$
with these linear combinations, i.e. these columns are now all-zero columns in a 3-subspace which replaces $X$.
These three linear combinations are performed on all the 3-subspaces of $\dS$ and replace the related columns
in all the 381 3-subspaces of $\dS$. By abuse of notation we call the new system also $\dS$. We note that after this was done, $Z_1$ was not
affected and it remains a 3-subspace of $\dS$. Hence, the two 3-subspaces
(starting with four all-zero columns and ending with such four columns, i.e. $Z_1$ and $Z_2$) can be
forced to be in $\dS$ which we do. As a consequence, all the consequences that we will derive
concerning the first $\ell$, $1 \leq \ell \leq 6$, columns of the 3-subspaces in $\dS$,
are also correct consequences concerning the last $\ell$ columns of these 3-subspaces.
Let $\T$ be the system formed from $\dS$ by performing puncturing three times
on the first three columns of all the 3-subspaces of $\dS$
(note that $\T$ is isomorphic to $\dS_2(2,3,7;4)$, but such system
was defined before only when the last columns are punctured).

Each pair of vectors from $\F_2^7$, one which starts with four \emph{zeroes} and one which ends with four
\emph{zeroes} appears together in exactly one 3-subspace of $\dS$. There are no three such vectors
in the same 3-subspace of $\dS$ since two such vectors
(with either four leading \emph{zeroes} or four \emph{zeroes} at the tail)
will sum to another such vector and the result will be a 1-subspace in either $\dS_2(2,3,7,;4)$ or $\T$, a contradiction.
Therefore, there are exactly $7 \cdot 7=49$ 3-subspaces with one vector who has
four leading \emph{zeroes} and one vector who has four \emph{zeroes} at the tail.
Let $\A$ be the set of 3-subspaces of $\dS$ which form the 140 2-subspaces in $\dS_2(2,3,7;4)$ and
let $\dB$ be the set of 3-subspaces of $\dS$ which form the 140 2-subspaces in $\T$.
Clearly,
$$
|\A|=|\dB|=140,~~~ |\A \cap \dB|=49,~~~|\A \setminus\dB|=|\dB \setminus \A|=91~~.
$$

Therefore, there are $381-(91+91+49+1+1)=148$ 3-subspaces in $\dS$ in which the projection
on the first four columns and the projection on the last four columns yield a 3-subspace
of $\F_2^4$.

Can we have another 3-subspace in $\dS$ with four all-zero columns. We cannot give a definite answer to this
question. In such a 3-subspace, two all-zero columns must be in the first three columns and the
two other all-zero columns in the last three columns,
which even make it possible for three more different 3-subspaces like this in $\dS$ for a total
of five 3-subspaces with four all-zero columns. A short discussion on two subspaces in $\dS$
with four all-zero columns which share two coordinates of all-zero columns will be given for any
prime power $q$ in the next subsection.

Based on the forced structure described in this section, one can start a computer search either to construct
the $q$-Fano plane for $q=2$. The outcome of such search is of great interest.

\subsection{The Structure of $\dS_q(2,3,7)$, General $q$}
\label{sec:q>1}

Throughout our discussion, let $\dS$ be a $q$-Steiner system $\dS_q(2,3,7)$.
We start with a uniform solution for the 3-punctured $q$-Steiner system
$\dS_q(2,3,7;4)$. Such a uniform solution, given in~\cite{Etz15} implies that $X_{0,4}=1$ which implies
that $X_{1,4}=0$, $X_{2,4}=q^2$, and $X_{3,4}=q^4(q-1)$. We set $X_{0,4}=1$ w.l.o.g.,
since in any system one subspace can be chosen w.l.o.g.~. Furthermore, in this case
$X_{0,4}=1$ implies that the design is uniform~\cite{Etz15}. $X_{0,4}=1$ implies that the 3-subspace whose
first four columns are all-zero columns is contained in $\dS$. Let $Z_1$ denote this
3-subspace in which the first four columns are all-zero. For symmetry let $Z_2$ denote
the 3-subspace of $\F_q^7$ in which the last four columns are all-zero columns.

In $\dS$, each 1-subspace of $\F_q^7$ contained in exactly $\frac{q^6-1}{q-1}$ 3-subspaces.
By puncturing the last coordinate of each 3-subspace of $\dS$,
all the 3-subspaces which contain the vector 0000001 will be punctured into
a spread, i.e. a $q$-Steiner system $\dS_q(1,2,6)$.

Next, in our exposition we exclude $Z_1$ from $\dS$ for the current paragraph.
Each 3-subspace of $\dS$ which contains a nonzero vector which starts with four \emph{zeroes}
is 3-punctured into a 2-subspace of $\F_q^4$. There are $\frac{q^3-1}{q-1}$ 1-subspaces which contain such vectors, each one is contained
in $q^2 (q^2+1)$ distinct 3-subspaces of $\dS \setminus \{ Z_1 \}$ (since the only 3-subspace of $\dS$ which
contains two such 1-subspaces is $Z_1$; thus, clearly two of them cannot be contained together in the same 3-subspace of $\dS \setminus \{ Z_1 \}$)
for a total of $q^2 (q^2+1) (q^2+q+1)$ such 3-subspaces which are punctured
into the $q^2 (q^2+1) (q^2+q+1)$ (non-distinct) 2-subspaces
of the 3-punctured $q$-Steiner system $\dS_q(2,3,7;4)$ derived from $\dS$.
Since each 1-subspace of $\F_q^7$ is contained in exactly one 3-subspace with each
other 1-subspace of $\F_q^7$, it follows that each such 1-subspace is responsible
for exactly $q^2 (q^2+1)$ 2-subspaces (some of them are identical) of the 3-punctured $q$-Steiner system $\dS_q(2,3,7;4)$.
Hence, each nonzero prefix of length 4 of vectors of $\F_q^7$ (there are $q^4-1$ such nonzero prefixes) appears $\frac{q^2 (q^2+1) \cdot (q^2-1)}{q^4-1}=q^2$
times in these $q^2 (q^2+1)$ 2-subspaces. There are many possible partitions to obtain such $q^2+q+1$ 3-punctured sets
of $q^2 (q^2+1)$ 2-subspaces for this purpose.
One suggestion was given in the recursive construction of~\cite{Etz15}.
It will be discussed again in Section~\ref{sec:2punctured}.

We continue by imposing w.l.o.g. a certain structure on $\dS$.
Many structures can be imposed, each one will imply a different specifications for $\dS$. The one we suggest
now seems to be one of the most successful ones. We already forced the 3-subspace $Z_1$
to be a 3-subspace in $\dS$. Now, we assume
that w.l.o.g. also the 3-subspace whose last four columns are all-zero, i.e. $Z_2$, is also a 3-subspace in $\dS$.

Next, we will show why it can be assumed that $Z_2 \in \dS$ without loss of generality.
In the 3-punctured $q$-Steiner system $\dS_q(2,3,7;4)$ each 3-subspace of $\F_q^4$ appears exactly $q^4 (q-1)$ times.
There are $\frac{q^4-1}{q-1}$ such 3-subspaces of $\F_q^4$, for one of them, the first three columns
form the unique 3-subspace in $\F_q^3$, and these three columns can be followed by the all-zero column (which is one
of the $q^3$ possible extensions of the unique 3-subspace of $\F_q^3$ by Lemma~\ref{lem:addNone}). Let $X$
be a 3-subspace of $\dS$ which is 3-punctured to these four columns. To continue
our discussion, we use Lemma~\ref{lem:linear_comb}.
Since the first three columns of $X$ have rank three it follows that we can
form some specific three linear combinations, containing the 5th, the 6th, and the 7th column of $X$, respectively.
Each such linear combination will sum
to \emph{zero} for the related column of $X$.
We replace the 5th, 6th, and 7th columns of $X$
with these linear combinations, i.e. these columns are now all-zero columns in a 3-subspace which replaces $X$.
These three linear combinations are performed and replace the related columns
in all the $(q^6+q^5+q^4+q^3+q^2+q+1)(q^2-q+1)$ 3-subspaces of $\dS$.
By abuse of notation we call the new system also $\dS$. We note that after this was done, $Z_1$ was not
affected and it remains a 3-subspace of $\dS$. Hence, the two 3-subspaces
(starting with four all-zero columns and ending with such four columns, i.e. $Z_1$ and $Z_2$) can be
forced to be in $\dS$ which we do. As a consequence, all the consequences that we will derive
concerning the first $\ell$, $1 \leq \ell \leq 6$, columns of the 3-subspaces in $\dS$,
are also correct consequences concerning the last $\ell$ columns of these 3-subspaces.
Let $\T$ be the system formed from $\dS$ by performing puncturing three times
on the first three columns of all the 3-subspaces of $\dS$ (note that $\T$ is isomorphic to $\dS_q(2,3,7;4)$, but such system
was defined before only when the last columns are punctured).

Each pair of 1-subspaces of $\F_q^7$ which contain vectors which start with four \emph{zeroes} and vectors which end with four
\emph{zeroes} appears together in exactly one 3-subspace of $\dS$. There are no three such linearly independent vectors
in the same 3-subspace of $\dS$ since two such vectors (with either four leading \emph{zeroes} or four \emph{zeroes} at the tail)
will sum to another such vector and the result will be a 1-subspace in either $\dS_q(2,3,7,;4)$ or $\T$, a contradiction.
Therefore, there are exactly $\frac{q^3-1}{q-1} \cdot \frac{q^3-1}{q-1}=(q^2+q+1)^2$ 3-subspaces which contain 1-subspace who has
a vector with four leading \emph{zeroes} and one 1-subspace which contains a vector who has four \emph{zeroes} at the tail.
Let $\A$ be the set of 3-subspaces of $\dS$ which form the $q^2 (q^2+1)(q^2+q+1)$ 2-subspaces in $\dS_q(2,3,7;4)$ and
let $\dB$ be the set of 3-subspaces of $\dS$ which form the $q^2 (q^2+1)(q^2+q+1)$ 2-subspaces in $\T$.
Clearly,
$$
|\A|=|\dB|=q^2 (q^2+1)(q^2+q+1),~~~ |\A \cap \dB|=(q^2+q+1)^2,~~~|\A \setminus\dB|=|\dB \setminus \A|=(q^2+q+1)(q^4-q-1)~~.
$$

Therefore, there are $(q^6+q^5+q^4+q^3+q^2+q+1)(q^2-q+1)-(2 \cdot (q^2+q+1)(q^4-q-1) + (q^2+q+1)^2+1+1)=q(q^7-q^5-q^4-2q^3+q^2+2q+2)$
3-subspaces in $\dS$ in which the projection
on the first four columns and the projection on the last four columns yield 3-subspace
in $\F_q^4$.

Finally, as we mentioned before, there are other possible 3-subspaces that can be imposed on $\dS$
in addition to $Z_1$. We will briefly mention one more such option. We are mainly interested in
3-subspaces which have four all-zero columns since we know the structure of the related
design formed by puncturing the three other columns. We claim that we can impose on $\dS$ to have
two such 3-subspaces (one of them is $Z_1$) and one which has two joint all-zero columns with $Z_1$
($Z_2$ has only one such joint column). Say such 3-subspace is $Z_3$ and it has all-zero columns in columns 1, 2, 5, and 6.
The proof that we can force $Z_1$ and $Z_3$ to be together 3-subspaces in $\dS$
is very similar to the one which forced $Z_1$ and $Z_2$ to be together 3-subspaces in $\dS$.
We consider the unique 2-subspace of $\F_q^4$ in which the first and the second columns
are the all-zero column. There are four 3-subspaces extended from this subspace. Let $Y$
be one of these four subspaces. We will prove that w.l.o.g. we can assume that in $Y$ the first vector from $\F_q^7$
starts with six \emph{zeroes}. The first four entries are \emph{zeroes} since it was extended from a 2-subspace.
If the last entry is \emph{zero} we can exchange the 7th column with either the 5th or the 6th column
(one which has nonzero entry in this vector). Now, if the 5th or the 6th column does not have
a \emph{zero} in the first vector we replace it with a linear combination of the related column with
the 7th column making it \emph{zero} as suggested in Lemma~\ref{lem:linear_comb}.
Now, we can use the 3rd and 4th column by using Lemma~\ref{lem:linear_comb} again to make
the 5th and 6th columns of $Y$ all-zero columns. All these linear combinations done on all
the 3-subspaces of $\dS$ do not affect $Z_1$ and the final result is that $Z_1$ and $Z_3$ are in $\dS$.
There are other such 3-subspaces that can be forced to be in $\dS$ together with $Z_1$ and they might
also help either to construct a $q$-Fano plane or to prove its nonexistence.

\section{The 2-punctured $q$-Steiner system $S_q(2,3,7;5)$}
\label{sec:2punctured}

In this section we present a possible structure for the $q$-Fano plane, namely,
we present a construction of a 2-punctured $q$-Steiner system $\dS_q(2,3,7;5)$.
We note that this is a possible substructure of the $q$-Steiner system $\dS_q(2,3,7)$
(first five columns of the system), but it is not forced like the syetems described in Section~\ref{sec:structure},
and hence, it might not be possible to complete the construction to the $q$-Fano plane, even
if the related $q$-Fano plane exists. The construction is based on extension of all the
subspaces of the 3-punctured $q$-Steiner system $\dS_q(2,3,7;4)$.

Let $\dS$ be a 3-punctured $q$-Steiner system $\dS_q(2,3,7;4)$ with the uniform solution,
found in~\cite{Etz15}, i.e. $X_{0,4}=1$, $X_{1,4}=0$, $X_{2,4}=q^2$ and $X_{3,4}=q^4(q-1)$.

Each 3-subspace of $\F_q^4$ in $\dS$ can be extended in $q^3$ different ways (see Lemma~\ref{lem:addNone}).
We use each such extension $q(q-1)$ times in $S_q(2,3,7;5)$, i.e. each one of the
$\sbinomq{4}{3} q^3$ such 3-subspaces of $\F_q^5$ will appear $q(q-1)$ times in our
constructed $\dS_q(2,3,7;5)$.

There are $q^2 (q^2+1)(q^2+q+1)$ 2-subspaces in $S_q(2,3,7;4)$, i.e. $(q^2+1)(q^2+q+1)$
distinct 2-subspaces for which each one appears $q^2$ times in $\dS_q(2,3,7;4)$.
From this set of 2-subspaces there are $q^2 (q^2+1)q^2$ 2-subspaces which will be extended
to 3-subspaces and $q^2 (q^2+1) (q+1)$ which will be extended to 2-subspaces.
There are a total of $(q^2+1)(q^2+q+1)$ distinct 2-subspaces in $\F_q^4$ that can be partitioned
into $q^2+q+1$ disjoint spreads, each one of size $q^2+1$~\cite{Bak76,Beu74}.
We partition these disjoint spreads into two sets, one set $A$ will
contain $q^2$ spreads and a second set $B$ will contain $q+1$ spreads.

Each 2-subspace in $\dS$ which is contained in a spread from the
set $A$ is extended in a unique way (see Lemma~\ref{lem:addOne}) to a 3-subspace of $S_q(2,3,7;5)$.
Thus, each such 3-subspace of $\F_q^5$ (there are $q^2 (q^2+1)$ such
3-subspaces) will appear $q^2$ times in our constructed $S_q(2,3,7;5)$.

Each 2-subspace in $\dS$ which is contained in a spread from the
set $B$ is extended in a $q^2$ ways (see Lemma~\ref{lem:addNone}) to 2-subspaces of $S_q(2,3,7;5)$.
Thus, each one of these $q^2$ new 2-subspaces of $\F_q^5$ (there are $(q+1)(q^2+1)q^2$
such 2-subspaces) will appear exactly once in our constructed $S_q(2,3,7;5)$.

\section{Conclusion}
\label{sec:conclude}

We have presented a framework which can help to determine whether a $q$-Fano
plane, i.e. a $q$-Steiner system $\dS_q(2,3,7)$ exists. For small values of $q$,
probably only for $q=2$ this might help to determine the existence of such system
using computer search. A short step, rather than a complete solution to the problem,
to make a progress in solving the existence problem, can be done in one of
the following directions:

\begin{enumerate}
\item Find a punctured $q$-Steiner system $\dS_q(2,3,7;6)$. First step in this
direction would be to consider $q=2$.

\item The subspaces of the set $\A \cup \dB$ might be a key for the whole construction
of $\dS_q (2,3,7)$. A first step might be to find the 231 subspaces of this set for $q=2$
and to extend this set of 231 3-subspaces with as many as possible 3-subspaces, such that no 2-subspace of $\F_2^7$ appears in more than one
3-subspace. Of course, the two subspaces with four all-zero columns in the first or last columns must be part of this construction.

\item Another small step forward will be to settle the possibility of more than two
3-subspaces with four all-zero columns. It is either to prove that no more than two 3-subspaces
with four all-zero columns cannot exist in $\dS_q(2,3,7)$ or to prove that we can assume the existence
of three such 3-subspaces without loss of generality. We note that there is always the possibility that such a system
can have three 3-subspaces, but not without loss of generality.
\end{enumerate}


\end{document}